\documentclass{fundam}

\usepackage{amsmath,amscd,amssymb,stmaryrd}

\usepackage{comment}
\usepackage{tikz}
\usetikzlibrary{positioning,arrows,calc}
\usepackage{xspace}
\usepgflibrary{shadings}
\usepackage{algorithm}
\usepackage{algpseudocode}

\usepackage{subcaption}
\usepackage{graphicx}

\usepackage{url}

\newcommand{\R}{\mathbb{R}}
\newcommand{\Q}{\mathbb{Q}}
\newcommand{\Z}{\mathbb{Z}}

\newcommand{\N}{\mathbb{N}}

\newcommand{\hex}{\mathrm{hex}}

\newcommand{\expectation}{\mathbb{E}}

\usepackage{hyperref}

\begin{document}

 \setcounter{page}{331}
 \publyear{24}
 \papernumber{2186}
  \volume{191}
  \issue{3-4}

\finalVersionForARXIV


\title{Finding Codes on Infinite Grids Automatically}

\author{Ville Salo\thanks{Address for correspondence: Department of Mathematics and Statistics,
                                          University of Turku, Finland}, \  Ilkka T\"orm\"a\thanks{Supported by the Academy of Finland under grant 346566.}
 \\
Department of Mathematics and Statistics \\
University of Turku, Finland \\
vosalo@utu.fi,  iatorm@utu.fi
}

\maketitle

\runninghead{V.\ Salo and  I.\ T\"orm\"a}{Finding Codes on Infinite Grids Automatically}

\vspace*{4mm}
\begin{abstract}
We apply automata theory and Karp's minimum mean weight cycle algorithm to minimum density problems in coding theory. Using this method, we find the new upper bound $53/126 \approx 0.4206$ for the minimum density of an identifying code on the infinite hexagonal grid, down from the previous record of $3/7 \approx 0.4286$.
\end{abstract}


\section{Introduction}

It is well-known in symbolic dynamics that for a subshift of finite type which has sufficient gluing properties (such as strong irreducibility) the minimal (or maximal) frequency of any letter in a configuration is uniformly computable given the description of the subshift. This observation in particular applies to many subshifts studied in the field of coding theory, such as variants of identifying codes and locating-dominating codes. Indeed, the minimal frequency is lower semicomputable for any subshift of finite type, and under sufficient gluing properties it is upper semicomputable, and often even approximated arbitrarily well by periodic points.

\medskip
This paper is dedicated to Iiro Honkala on his 60th birthday. As Iiro has worked on such density problems, and they are really problems about symbolic dynamics (and ergodic theory) in thin disguise, this felt like a natural common point of interest to look into. In particular, Ville Junnila suggested to us the study of the minimal density of identifying codes on the so-called hexagonal grid. We will refer to this density as $\alpha_0$ in this article. The lowest-density code known so far was found by Cohen, Honkala, Lobstein and Zemor in 2000 \cite{CoHoLoZe00}, with density $3/7 \approx 0.4286$ (showing $\alpha_0 \leq 3/7$).

\medskip
The best known lower bound $\alpha_0 \geq 23/55 \approx 0.41818$ was found in 2014 by Derrick Stolee, who automated the so-called discharging method through linear programming \cite{St14}. Code construction has also previously been attacked with computer methods, in particular \cite{ChHuLo02} used a heuristic method to find codes for many identifying code problems in regular gridlike graphs.
However, this method was not able to improve on the upper bound $3/7$ for the hexagonal grid.

\begin{figure}[htp]
\begin{center}
\begin{tikzpicture}[scale=0.34]
\draw[fill=black!30!white] (9.5,-0.4) rectangle (35.5,1.1);
\draw[fill=black!30!white] (-0.5,1.1) rectangle (25.5,2.5);
\draw[fill=black!30!white] (15.5,2.5) rectangle (41.5,3.9);
\draw[fill=black!30!white] (5.5,3.9) rectangle (31.5,5.3);
\input{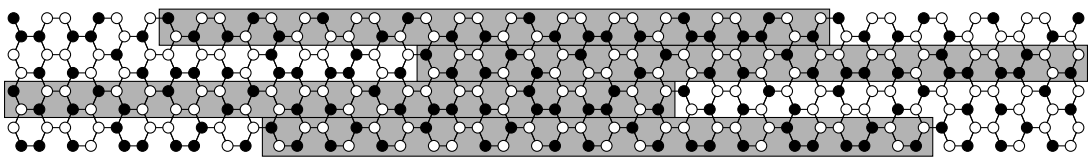}
\end{tikzpicture}
\end{center}
\caption{An identifying code with density $11/26$ on the hexagonal grid, found using the SAT solver Gluecard.}
\label{fig:sharpSAT}
\end{figure}

\medskip
We tried two new methods for lowering $\alpha_0$. First, we naively input this problem into the SAT-solver Gluecard, accessed through the PySAT library \cite{IgMoMa18}, to directly compute the minimal densities of periodic configurations of small periods. This was already successful: we found the code of density $11/26 \approx 0.4231$ shown in Figure~\ref{fig:sharpSAT}, reducing the gap between Stolee's lower bound and the known upper bound from about $0.0103$ to about $0.0049$.

\medskip
Second, as the configurations with one fixed period $\vec v$ form in essence a $\Z$-subshift of finite type, one can compute the exact minimal density of a configuration with period $\vec v$ using simple graph/automata theory, as it is always reached by a cycle. Specifically, we constructed a finite-state automaton for the $\vec v$-periodic configurations for various $\vec v$, recognizing the ``strip'' between periods, and used Karp's minimum mean weight cycle algorithm \cite{Ka78} to compute the exact minimal density. We ran the program on the Dione cluster of the FGCI (Finnish Grid and Cloud Infrastructure) project, which has 80 CPUs and about 200 gigabytes of RAM per computation node. The Python source code is available at \cite{githubhex}.

\def\eps{0.002}
\def\kkk{0}
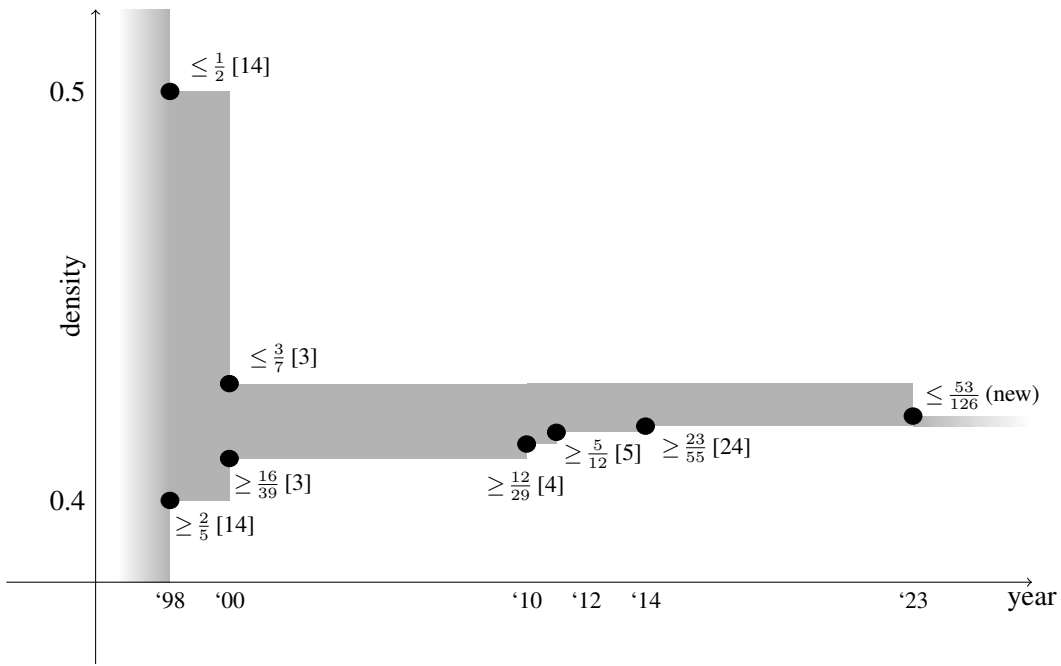
\begin{figure}[!ht]
\begin{center}
\scalebox{0.98}{
\begin{tikzpicture}[scale=50,yscale=1.1,xscale=1.2]

\fill[left color=white, right color=black!30!white, middle color=black!15!white] (-0.025,0.38) rectangle (-0.2/15,0.52);

\draw (-0.03,0.36) edge[->] (-0.03,0.52);
\draw (-0.05,0.38) edge[->] (2.7/15,0.38);
\node[left] () at (-0.03,0.4) {0.4};
\node[left] () at (-0.03,0.5) {0.5};
\node[left] () at (-0.03,0.45) {\rotatebox{90}{density}};
\node[below] () at (2.7/15,0.38) {year};

\node[below] () at (-0.2/15,0.38) {\footnotesize{`98}};
\node[below] () at (-0/15,0.38) {\footnotesize{`00}};
\node[below] () at (1/15,0.38) {\footnotesize{`10}};
\node[below] () at (1.2/15,0.38) {\footnotesize{`12}};
\node[below] () at (1.4/15,0.38) {\footnotesize{`14}};
\node[below] () at (2.3/15,0.38) {\footnotesize{`23}};

\node[label={[label distance=-0.3cm]below right:{\footnotesize$\geq\!\frac25$\,\cite{KaChLe98}}}] (a) at (-0.2/15, 2/5) {}; 
\node[label={[label distance=\kkk]above right:{\footnotesize$\leq\!\frac12$\,\cite{KaChLe98}}}] (b) at (-0.2/15, 1/2) {}; 
\node[label={[label distance=-0.3cm]below right:{\footnotesize$\geq\!\frac{16}{39}$\,\cite{CoHoLoZe00}}}] (c) at (0, 16/39) {}; 
\node[label={[label distance=\kkk]above right:{\footnotesize$\leq\!\frac37$\,\cite{CoHoLoZe00}}}] (d) at (0, 3/7) {}; 
\node[label={[label distance=0]below:\rotatebox{0}{\footnotesize$\geq\!\frac{12}{29}$\,\cite{CrYu10}}}] (e) at (1/15, 12/29) {}; 
\node[label={[label distance=-11]below right:\rotatebox{0}{\footnotesize$\geq\!\frac{5}{12}$\,\cite{CuYu13}}}] (f) at (1.13/15, 5/12) {}; 
\node[label={[label distance=-11]below right:\rotatebox{0}{\footnotesize$\geq\!\frac{23}{55}$\,\cite{St14}}}] (g) at (1.46/15, 23/55) {}; 
\node[label={[label distance=-0.1cm]above right:{\footnotesize$\leq\!\frac{53}{126}$\,(new)}}] (h) at (2.3/15, 53/126) {}; 

\fill[black!30!white] (-0.2/15, 2/5) rectangle (0,1/2);
\fill[black!30!white] (-0.2/15, 16/39) rectangle (1/15,3/7);
\fill[black!30!white] (1/15, 12/29) rectangle (1.1/15,3/7);
\fill[black!30!white] (1.1/15, 5/12) rectangle (1.4/15,3/7);
\fill[black!30!white] (1.4/15, 23/55) rectangle (2.3/15,3/7);

\fill[left color=black!30!white, right color=white, middle color=black!15!white] (2.3/15,23/55) rectangle (2.7/15,53/126);

\filldraw (-0.2/15, 2/5) circle (0.002);
\filldraw (-0.2/15, 1/2) circle (0.002);
\filldraw (0, 16/39) circle (0.002);
\filldraw (0, 3/7) circle (0.002);
\filldraw (1/15, 12/29) circle (0.002);
\filldraw (1.1/15, 5/12) circle (0.002);
\filldraw (1.4/15, 23/55) circle (0.002);
\filldraw (2.3/15, 53/126) circle (0.002);
\end{tikzpicture}  }
\end{center}\vspace*{-4mm}
\caption{Evolution of the known lower and upper bounds on the minimal density of an identifying code on the hex grid. Years are based on publication years for \cite{KaChLe98} and \cite{CoHoLoZe00}, and arXiv years for others. The problem was introduced in \cite{KaChLe98}.}\label{fig:Years}
\end{figure}

\medskip
This method found a better code than the SAT solution, with density $53/126 \approx 0.4206$, reducing the gap between Stolee's lower bound and the known upper bound to about $0.0025$. The evolution of the bounds is illustrated in Figure~\ref{fig:Years}. The automata method turned out superior to the SAT solver (at least on this one problem), as the solver was not able to reproduce the optimal solution even given the periods.

\begin{theorem}
\label{thm:NewRecord}
The minimal density identifying code on the hexagonal grid is at most $53/126$.
\end{theorem}

\begin{figure}[htp]
\begin{center}
\includegraphics[width=12.5cm]{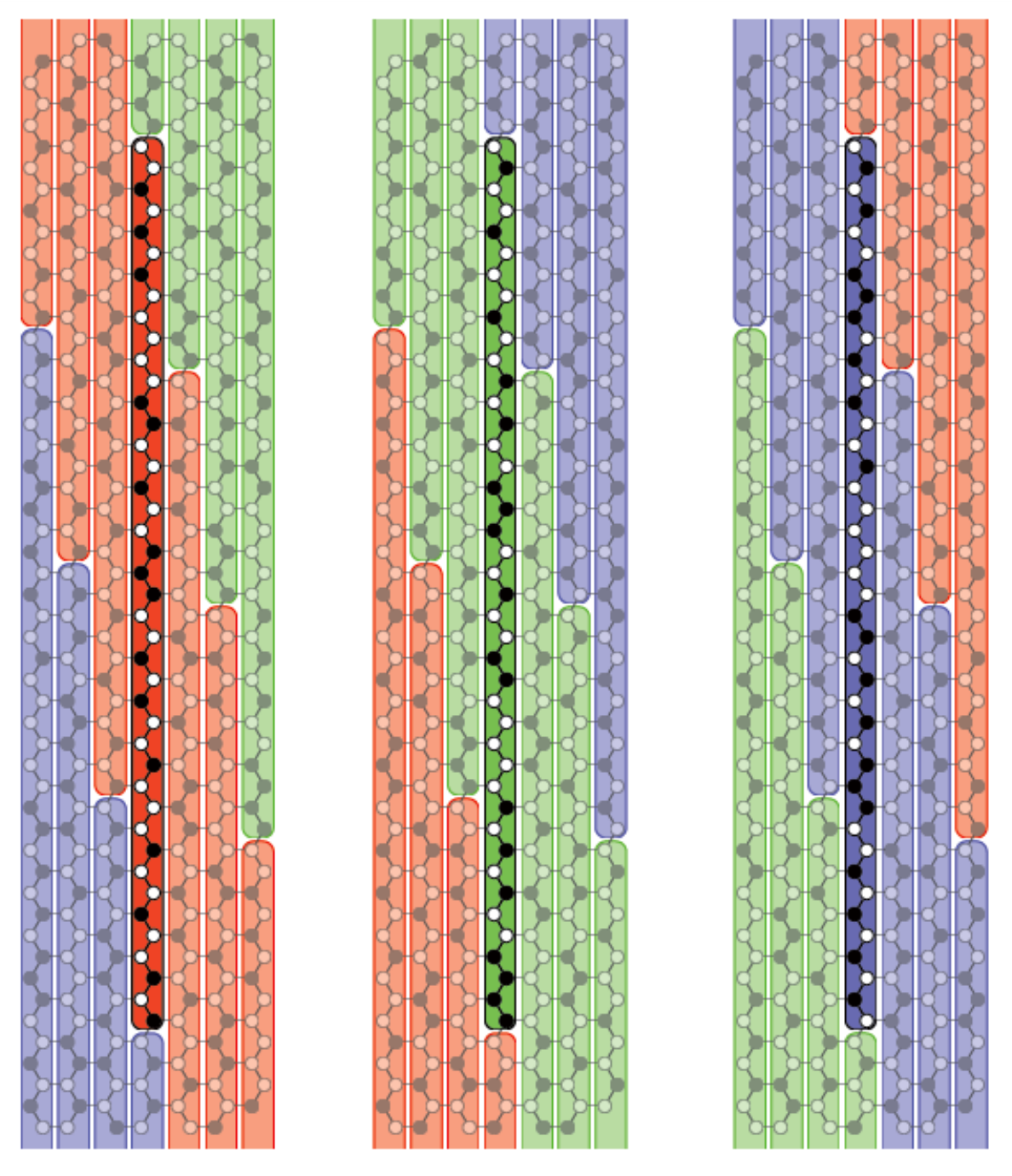}
\end{center}
\caption{An identifying code with density $53/126$ on the hexagonal grid found using automata theory and Karp's algorithm.}
\label{fig:NewRecord}
\end{figure}

\begin{proof}
In Figure~\ref{fig:NewRecord}, we exhibit an identifying code with periods $(0,126)$ and $(1,-11)$.
\end{proof}

We also looked at four other problems from the literature with the automata-theoretic method, namely
\begin{itemize}
\itemsep=1.2pt
\item the minimum density of a radius-$2$ locating-dominating code on the king grid, $\alpha_1$; previously known bounds $\alpha_1 \in [1/10, 1/8]$ from \cite{ChHuLo02,Pe13,Pe16},
\item the minimum density of a redundant radius-1 locating-dominating code on the king grid, $\alpha_2$; previously known bounds $\alpha_2 \in [3/11, 5/16]$ from \cite{DeSu22},
   \eject
\item the minimum density of a radius-$2$ identifying code for the square grid, $\alpha_3$; previously known bounds $\alpha_3 \in [6/35, 5/29]$ from \cite{Ju13,HoLo02}, and
\item the minimum density of a radius $2$ identifying code for the triangular grid, $\alpha_4$; previously known bounds $\alpha_4 \in [2/15, 1/6]$ from \cite{ChHoHuLo01}.
\end{itemize}
In each case, we ran the search for all periods that the cluster could handle before running out of memory or time. We were not able to improve on any of the densities, which could mean that our method does not apply well to these codes, or alternatively might suggest that these codes are (close to) optimal and more work is needed on the side of lower bounds.

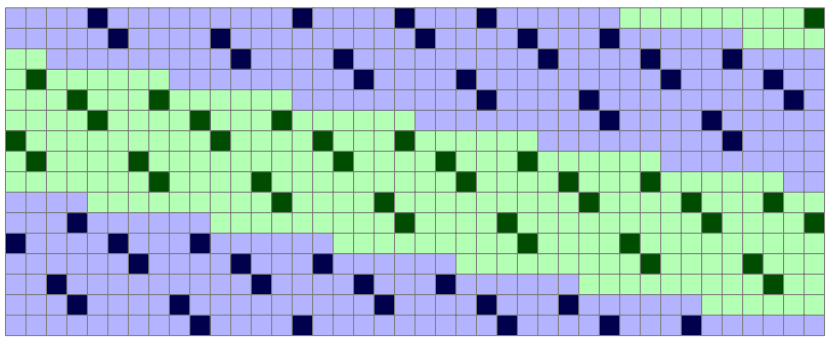
\begin{figure}[htp]
\begin{center}
\begin{tikzpicture}[scale=0.27]
\fill (9, 0) rectangle (10, 1);\fill (14, 0) rectangle (15, 1);\fill (24, 0) rectangle (25, 1);\fill (29, 0) rectangle (30, 1);\fill (33, 0) rectangle (34, 1);\fill (3, 1) rectangle (4, 2);\fill (8, 1) rectangle (9, 2);\fill (18, 1) rectangle (19, 2);\fill (23, 1) rectangle (24, 2);\fill (27, 1) rectangle (28, 2);\fill (2, 2) rectangle (3, 3);\fill (12, 2) rectangle (13, 3);\fill (17, 2) rectangle (18, 3);\fill (21, 2) rectangle (22, 3);\fill (37, 2) rectangle (38, 3);\fill (6, 3) rectangle (7, 4);\fill (11, 3) rectangle (12, 4);\fill (15, 3) rectangle (16, 4);\fill (31, 3) rectangle (32, 4);\fill (36, 3) rectangle (37, 4);\fill (0, 4) rectangle (1, 5);\fill (5, 4) rectangle (6, 5);\fill (9, 4) rectangle (10, 5);\fill (25, 4) rectangle (26, 5);\fill (30, 4) rectangle (31, 5);\fill (3, 5) rectangle (4, 6);\fill (19, 5) rectangle (20, 6);\fill (24, 5) rectangle (25, 6);\fill (34, 5) rectangle (35, 6);\fill (39, 5) rectangle (40, 6);\fill (13, 6) rectangle (14, 7);\fill (18, 6) rectangle (19, 7);\fill (28, 6) rectangle (29, 7);\fill (33, 6) rectangle (34, 7);\fill (37, 6) rectangle (38, 7);\fill (7, 7) rectangle (8, 8);\fill (12, 7) rectangle (13, 8);\fill (22, 7) rectangle (23, 8);\fill (27, 7) rectangle (28, 8);\fill (31, 7) rectangle (32, 8);\fill (1, 8) rectangle (2, 9);\fill (6, 8) rectangle (7, 9);\fill (16, 8) rectangle (17, 9);\fill (21, 8) rectangle (22, 9);\fill (25, 8) rectangle (26, 9);\fill (0, 9) rectangle (1, 10);\fill (10, 9) rectangle (11, 10);\fill (15, 9) rectangle (16, 10);\fill (19, 9) rectangle (20, 10);\fill (35, 9) rectangle (36, 10);\fill (4, 10) rectangle (5, 11);\fill (9, 10) rectangle (10, 11);\fill (13, 10) rectangle (14, 11);\fill (29, 10) rectangle (30, 11);\fill (34, 10) rectangle (35, 11);\fill (3, 11) rectangle (4, 12);\fill (7, 11) rectangle (8, 12);\fill (23, 11) rectangle (24, 12);\fill (28, 11) rectangle (29, 12);\fill (38, 11) rectangle (39, 12);\fill (1, 12) rectangle (2, 13);\fill (17, 12) rectangle (18, 13);\fill (22, 12) rectangle (23, 13);\fill (32, 12) rectangle (33, 13);\fill (37, 12) rectangle (38, 13);\fill (11, 13) rectangle (12, 14);\fill (16, 13) rectangle (17, 14);\fill (26, 13) rectangle (27, 14);\fill (31, 13) rectangle (32, 14);\fill (35, 13) rectangle (36, 14);\fill (5, 14) rectangle (6, 15);\fill (10, 14) rectangle (11, 15);\fill (20, 14) rectangle (21, 15);\fill (25, 14) rectangle (26, 15);\fill (29, 14) rectangle (30, 15);\fill (4, 15) rectangle (5, 16);\fill (14, 15) rectangle (15, 16);\fill (19, 15) rectangle (20, 16);\fill (23, 15) rectangle (24, 16);\fill (39, 15) rectangle (40, 16);\fill[opacity=0.3,color=blue] (0,0) rectangle (40,0+1);\fill[opacity=0.3,color=blue] (0,1) rectangle (34,1+1);\fill[opacity=0.3,color=green] (34,1) rectangle (40,1+1);\fill[opacity=0.3,color=blue] (0,2) rectangle (28,2+1);\fill[opacity=0.3,color=green] (28,2) rectangle (40,2+1);\fill[opacity=0.3,color=blue] (0,3) rectangle (22,3+1);\fill[opacity=0.3,color=green] (22,3) rectangle (40,3+1);\fill[opacity=0.3,color=blue] (0,4) rectangle (16,4+1);\fill[opacity=0.3,color=green] (16,4) rectangle (40,4+1);\fill[opacity=0.3,color=blue] (0,5) rectangle (10,5+1);\fill[opacity=0.3,color=green] (10,5) rectangle (40,5+1);\fill[opacity=0.3,color=blue] (0,6) rectangle (4,6+1);\fill[opacity=0.3,color=green] (4,6) rectangle (40,6+1);\fill[opacity=0.3,color=green] (0,7) rectangle (38,7+1);\fill[opacity=0.3,color=blue] (38,7) rectangle (40,7+1);\fill[opacity=0.3,color=green] (0,8) rectangle (32,8+1);\fill[opacity=0.3,color=blue] (32,8) rectangle (40,8+1);\fill[opacity=0.3,color=green] (0,9) rectangle (26,9+1);\fill[opacity=0.3,color=blue] (26,9) rectangle (40,9+1);\fill[opacity=0.3,color=green] (0,10) rectangle (20,10+1);\fill[opacity=0.3,color=blue] (20,10) rectangle (40,10+1);\fill[opacity=0.3,color=green] (0,11) rectangle (14,11+1);\fill[opacity=0.3,color=blue] (14,11) rectangle (40,11+1);\fill[opacity=0.3,color=green] (0,12) rectangle (8,12+1);\fill[opacity=0.3,color=blue] (8,12) rectangle (40,12+1);\fill[opacity=0.3,color=green] (0,13) rectangle (2,13+1);\fill[opacity=0.3,color=blue] (2,13) rectangle (40,13+1);\fill[opacity=0.3,color=blue] (0,14) rectangle (36,14+1);\fill[opacity=0.3,color=green] (36,14) rectangle (40,14+1);\fill[opacity=0.3,color=blue] (0,15) rectangle (30,15+1);\fill[opacity=0.3,color=green] (30,15) rectangle (40,15+1);\draw[gray] (0,0) grid (40,16);
\end{tikzpicture}
\end{center}
\caption{A new locating-dominating code of density $1/8$ for the radius-$2$ king grid. The colors green and blue are included to guide the eye and show the repeating pattern (in the form discovered by the algorithm, i.e.\ a row of width $40$).}
\label{fig:alpha1}
\end{figure}

\medskip
The most interesting discovery we made from these searches is the new locating-dominating code of density $1/8$ for $\alpha_1$, which has periods $(40,0)$ and $(-6,1)$, and took 87 hours of computation time on the Dione cluster to find. This code is pictured in Figure~\ref{fig:alpha1}. The previously known code from \cite{ChHuLo02} simply repeats the block
\begin{center}
\tikz[scale=0.3,baseline=0]{\draw[gray] (0,0) grid (4,4); \fill (0,0) rectangle (1,1); \fill (2,2) rectangle (3,3);}
\end{center}
with horizontal and vertical period $4$. Interestingly, the old code is also identifying, but our new code in Figure~\ref{fig:alpha1} is not. 

\section{Preliminaries}

We recall some definitions from multidimensional symbolic dynamics.
Let $A$ be a finite alphabet.
A $d$-dimensional \emph{pattern} is a mapping $P : D \to A$ for a finite or infinite \emph{domain} $D = D(P) \subset \Z^d$.
The $d$-dimensional \emph{full shift} over $A$ is the set $A^{\Z^d}$ of full patterns equipped with the prodiscrete topology, which is generated by the clopen \emph{cylinder sets} $[P] = \{ x \in A^{\Z^d} \mid x|_{D(P)} = P \}$ for all finite patterns $P$.
The group $\Z^d$ acts on $A^{\Z^d}$ by the \emph{shifts} $\sigma_{\vec v} : A^{\Z^d} \to A^{\Z^d}$, defined by $\sigma_{\vec v}(x)_{\vec n} = x_{\vec n + \vec v}$, which are homeomorphisms.
If $d = 1$, we denote $\sigma_1 = \sigma$ and call it the \emph{left shift}.
For $L \subset \Z^d$, a configuration $x \in A^{\Z^d}$ is \emph{$L$-periodic} if $x = \sigma_{\vec v}(x)$ for all $\vec v \in L$.
If $L$ is a rank-$d$ subgroup of $\Z^d$, we say $x$ is \emph{totally periodic}.

\medskip
A \emph{subshift} is a topologically closed $\sigma$-invariant subset $X \subset A^{\Z^d}$. The elements $x \in X$ are called \emph{configurations}. Every subshift is defined by some set of \emph{forbidden finite patterns} $F$:
\[
X = \bigcap_{P \in F} \bigcap_{\vec v \in \Z^d} \sigma_{\vec v}(A^{\Z^d} \setminus [P]).
\]
If $F$ can be chosen finite, then $X$ is a \emph{shift of finite type} (SFT). The finite patterns $P$ such that $[P] \cap X \neq \emptyset$ form the \emph{language} of the subshift $X$.

\medskip
A subshift $X \subset A^{\Z^d}$ is \emph{strongly irreducible} if there exists a \emph{gluing radius} $r \geq 0$ such that for all configurations $x, y \in X$ and all sets $U, V \subset \Z^d$ with $\min \{ \| \vec v - \vec w \| \mid \vec v \in U, \vec w \in V \} \geq r$, there exists a third configuration $z \in X$ with $z|_U = x|_U$ and $z|_V = y|_V$.
It has the \emph{uniform filling property} (or \emph{UFP}) if this holds whenever $U$ is a hypercube, and it is \emph{block gluing} if this holds whenever both $U$ and $V$ are hypercubes.

\medskip
We use several types of graphs.
If $V$ is a (countable, but possibly infinite) set and $E \subset \{(u, v) \mid u, v \in V, u \neq v\}$, then $G = (V, E)$ is a \emph{directed graph}. If $(u,v) \in E$ implies $(v,u) \in E$, we say $G$ is an \emph{undirected graph}. When an undirected graph is understood from context, $v \in V$ and $r \geq 1$, we write $B_r(v)$ for the closed ball of radius $r$ in $(V, E)$ with the path metric, or explicitly
\[ B_r(v) = \{u \in V \mid \exists n \leq r: \exists v = v_0, v_1, \ldots, v_n = u: \forall i: (v_i, v_{i+1}) \in E \}. \]
A \emph{weighted graph} is a triple $G = (V,E,\lambda)$ where $(V,E)$ is a graph and $\lambda : E \to \R$ associates a real weight to each edge.
In our applications, weighted graphs are always finite and directed, and they should be thought of as finite state automata.

\medskip
A \emph{radius-$r$ identifying code} \cite{KaChLe98} on an undirected graph $G = (V,E)$ is a subset $C \subset V$, whose elements are called \emph{codewords}, satisfying the following conditions:
\begin{itemize}
\item $B_r(u) \cap C \neq \emptyset$ for all $u \in V$,
\item $B_r(u) \cap C \neq B_r(v) \cap C$ for all distinct $u, v \in V$.
\end{itemize}
The set $B_r(u) \cap C$ is called the \emph{identifying set} of $u$ in $C$.
If we only require the second condition for distinct $u, v \in V \setminus C$, then the code is called a \emph{radius-$r$ locating-dominating code} \cite{Sl88}. By definition every identifying code is locating-dominating, but the converse is typically not true. When the radius is not mentioned, it is taken to be $1$. A \emph{redundant radius-$r$ locating-dominating code} is a radius-$r$ locating-dominating code that stays radius-$r$ locating-dominating even if any single codeword is removed. These were first studied in \cite{Sl02}.

Like in the case of subshifts, the set $2^V$ of all functions $x : V \to \{0, 1\}$ is considered with the topology with clopen basis given by cylinders $[P] = \{x \in 2^V \mid x|_D = P \}$ for patterns $P : D \to \{0,1\}$ with finite domain $D = D(P) \subset V$, and is homeomorphic to the Cantor space when $V$ is countably infinite. It is easy to see that radius-$r$ identifying codes are a closed subset of $2^V$, and the same is true for (redundant) locating-dominating codes.

We recall a few basic facts about abelian groups: every finitely-generated group $G$ which is torsion-free (meaning the identity element is not a proper power of another element) is isomorphic to $\Z^d$ under addition, for some $d$. In particular, all subgroups $L \leq \Z^d$ are isomorphic to $\Z^e$ for some $e \leq d$ ($e$ cannot be larger than $d$ by basic algebra, or even by considering growth rates). Since $\Z^d$ and $\Z^{d'}$ are not isomorphic when $d \neq d'$, the value of $e$ is uniquely determined by the subgroup $L$, and this is called the \emph{rank} of $L$. If $L \leq \Z^d$ is of rank $e < d$, then there exists $\vec v \in \Z^d$ such that $\langle L, \vec v \rangle$ has rank $e+1$. Indeed, one can pick any $\vec v$ outside the convex hull of $L$ (which is an $e$-dimensional subspace of $\R^d$). A rank-$d$ subgroup of $\Z^d$ is sometimes called a \emph{lattice}. If $X$ is a subshift, we say $x \in X$ is \emph{totally periodic} if the set of shifts that fix it form a lattice.

We also need basic definitions from ergodic theory. For general references on ergodic theory see e.g.\ \cite{DeGrSi76,Wa82,Ge11}. Here, we only consider ergodic theory on subshifts on groups $\Z^d$. For the purpose of measure theory, a subshift $X \subset A^{\Z^d}$ is always considered with its Borel $\sigma$-algebra (generated by the clopen sets). By a \emph{measure} we always mean a Borel probability measure $\mu$, i.e.\ $\mu(X) = 1$ and $\mu$ gives a measure precisely to the Borel sets. We write $\mathcal{M}$ for the set of all measures on $X$, and $\mathcal{M}_\sigma$ for the shift-invariant ones, i.e.\ those $\mu$ satisfying $\mu(\sigma_{\vec v}(B)) = \mu(B)$ for all Borel sets $B$ and all vectors $\vec v \in \Z^d$. The phrase ``for $\mu$-almost all $x$'' means ``for all $x \in X_0$ for some set $X_0$ with $\mu(X_0)=1$''. A measure is \emph{ergodic} if $\mu(B) \in \{0,1\}$ for all Borel sets $B \subset X$ satisfying $\sigma_{\vec v}(B) = B$.

If $X, Y$ are topological spaces, $f : X \to Y$ is a \emph{simple continuous function} if it is continuous and for some finite partition $C_1, \ldots, C_n$ of $X$ and values $y_1, \ldots, y_n \in Y$, we have $f|_{C_i} \equiv y_i$. Note that if $Y$ is Hausdorff, the $C_i$ must form a clopen partition. If $X$ is a $d$-dimensional subshift, a simple continuous function $f : X \to Y$ has a \emph{radius} $R$ such that $f(x)$ only depends on $x|_{B_R}$, where $B_R \subset \Z^d$ denotes the ball of radius $R$ with respect to Manhattan distance.

If $X$ is a subshift and $\mu$ is a measure on $X$, then a function $f : X \to \R_{\geq 0}$ is $L^1(\mu)$ if $f$ is measurable (preimages of intervals are Borel sets) and the Lebesgue integral $\int f d\mu$ is finite. It is $L^\infty(\mu)$ (or just $L^\infty$) if it is measurable and for some $X_0 \subset X$ with $\mu(X_0) = 1$ and $r < \infty$ we have $f(x) < r$ for all $x \in X_0$. All simple continuous functions to $\R_{\geq 0}$ are of course $L^1$ and $L^{\infty}$.

If $f : X \to \R_{\geq 0}$ is an $L^1(\mu)$ function and $\mathcal{I}$ is a $\sigma$-algebra contained in the Borel sets of $X$, then $\expectation_\mu(f | \mathcal{I})$ denotes the expectation of $f$ on $\mathcal{I}$, i.e.\ the unique $\mathcal{I}$-measurable (meaning preimages of intervals are in $\mathcal{I}$) function $g : X \to \R_{\geq 0}$ whose integral w.r.t.\ $\mu$ agrees with that of $f$ on every $B \in \mathcal{I}$. Such an expectation always exists, and is unique in the sense that any two expectations of $f$ on $\mathcal{I}$ agree in values on a set of measure $1$.

We also need a few notions from computability; for a general reference see e.g.\ \cite{Ro87}. A subshift is \emph{effectively closed}, if it can be defined by a computably enumerable set of forbidden patterns, i.e.\ a Turing machine can output some defining set of forbidden patterns (in any order). Equivalently, the entire complement of the language of the subshift can be enumerated. Such a machine is called an \emph{effective description (by forbidden patterns)} of the subshift. If a subshift is explicitly assumed to be an SFT, its description is a finite list of forbidden patterns. A real number $\alpha \in \R$ is \emph{lower (upper) semicomputable} if there exists an algorithm that outputs (in any order) all rational numbers strictly smaller (greater) than $\alpha$. If a real number is both lower and upper computable, we say it is computable.

Consider a set $L$ of finite words over some finite alphabet and a function $F : L \to \R$.
In our applications, $L$ will typically be a set of effective descriptions of some subclass of effectively closed subshifts, and $F(w)$ will be the minimum density of the subshift described by the Turing machine with encoding $w \in L$.
We say $F(w)$ is \emph{lower (upper) semicomputable uniformly in $w$} if there exists a single algorithm that, given $w$, outputs all rational numbers strictly smaller (greater) than $F(w)$.
If $F(w)$ is uniformly both lower and upper semicomputable in $w$, we say it is \emph{computable uniformly in $w$}.
If $F(w)$ is always rational and there exists an algorithm that, given $w$, computes it in finite time, we say $F(w)$ is \emph{exactly computable from $w$}.
Note that being exactly computable from $w$ is a strictly stronger property than being rational and computable uniformly in $w$.
In all of the above cases, the algorithm is allowed to do anything on input words that are not in $L$, even producing no output and/or failing to halt.

\section{Densities in subshifts}

\begin{definition}
\label{def:density}
Let $X \subset A^{\Z^d}$ be a subshift and $f : X \to \R$ a simple continuous function. The \emph{$f$-density} of $x \in X$ is
\[ \delta_f(x) = \lim_n \frac{\sum_{\vec v \in [-n,n]^d} f(\sigma_{\vec v}(x))}{(2n+1)^d} \]
when this limit exists.
The \emph{minimum density} of $f$ in $X$ is $\delta_f(X) = \inf_{x \in X} \delta_f(x)$, where for those $x \in X$ for which the $\delta_f(x)$ does not exist, we evaluate it to $\infty$.
\end{definition}

We immediately note that it cannot happen that the density is undefined for all points, i.e.\ the infimum does not take value $\infty$ at all points. This follows from Lemma~\ref{lem:DensCharacterization} below.

\medskip
For $a \in A$, the function
\[
 f(x) = \left\{\begin{array}{ll}
1 & \mbox{if } x_{\vec 0} = a \\
0 & \mbox{otherwise,}
\end{array}\right.
\]
is a simple continuous function, allowing us to define the density and minimum density of a single symbol.

\medskip
Suppose that $\Z^d$ acts freely on an undirected graph $G = (V,E)$ by graph automorphisms in such a way that the vertices have a finite number of distinct orbits.
We call such graphs \emph{gridlike}.
Let $D \subset V$ be a set of representatives for these orbits.
Then $V$ is the disjoint union $\bigcup_{\vec v \in \Z^d} (\vec v \cdot D)$, and the set $2^V$ of all functions $x : V \to \{0, 1\}$ is isomorphic to the full shift $(2^D)^{\Z^d}$ in the obvious way. 
Furthermore, the set of radius-$r$ identifying codes corresponds to a strongly irreducible subshift of finite type $X \subset (2^D)^{\Z^d}$, and the same is true for the other types of codes as well.
In this way, properties of SFTs can be transferred to those of codes on gridlike graphs; in particular, Definition \ref{def:density} gives us a notion of the density of 1s in a code $x : V \to \{0, 1\}$, and the minimum density of all codes on $G$.


\begin{remark}
In our definition of density, the hypercubes $[-n,n]^d$ could be replaced by any nested sequence of finite sets $N_1 \subset N_2 \subset \cdots \subset \Z^d$ with $\bigcup_i N_i = \Z^d$ and which satisfy the \emph{F\o{}lner condition}
\[
\lim_{i \to \infty} \frac{|N_i \triangle (N_i + \vec v)|}{|N_i|} = 0
\]
for all $\vec v \in \Z^d$, and this would lead to the same value for $\delta_f(X)$ and would not change any of the following results.
This includes the common definition via balls with respect to the path metric for gridlike graphs.
\end{remark}

We recall some basic results in ergodic theory (which hold in much greater generality). The use of these results is to replace various spatial averages by measures, which simplifies arguments where the geometric areas where density is computed can have different shapes.

The first result is a version of the famous pointwise ergodic theorem, and is a direct specialization of \cite[Theorem~14.A8]{Ge11}. For this theorem, consider any sequence of hypercubes $H_n = \vec v_n + [1, m_n]^d$ with $m_n \longrightarrow \infty$, such that $H_n \subset H_{n+1}$ for all $n$, and for $f : X \to \R_{\geq 0}$ define $(R_n f) (x) = \frac{1}{m_n^d}\sum_{\vec v \in H_n} f(\sigma_{\vec v}(x))$.

\begin{theorem}
\label{thm:ergodic}
Let $X \subset A^{\Z^d}$ be a subshift, and $\mathcal{I}$ the family of shift-invariant Borel sets. Then for any simple continuous function $f : X \to \R_{\geq 0}$ we have $\lim_n (R_n f)(x) = \expectation_\mu(f | \mathcal{I})(x)$ for $\mu$-almost all $x \in X$. 
\end{theorem}

\begin{corollary}
Suppose that $X$ is a $\Z^d$-subshift with invariant ergodic measure $\mu$, and $f : X \to \R_{\geq 0}$ is a simple continuous function. Letting $\alpha = \int f d\mu$, we have $(R_n f)(x) \longrightarrow \alpha$ for $\mu$-almost all $x \in X$.
\end{corollary}

\begin{proof}
By the previous theorem, $(R_n f)(x) \longrightarrow \expectation_\mu(f | \mathcal{I})(x)$ holds for $\mu$-almost all $x \in X$. The function $\expectation_\mu(f | \mathcal{I})$ is shift-invariant $\mu$-almost everywhere, so since $\mu$ is ergodic, it must be a constant function almost everywhere, i.e.\ $f(x) = \alpha$ for $\mu$-almost all $x$. We recall the argument: consider $A = f^{-1}(r)$. We have $\mu(\sigma_{\vec v}(A) \Delta A) = 0$ since $f$ is shift-invariant $\mu$-almost everywhere, where $\Delta$ denotes symmetric difference, so also $\mu((\bigcup_{\vec v \in \Z^d}\sigma_{\vec v}(A)) \Delta A) = 0$. Since $\bigcup_{\vec v \in \Z^d}(\sigma_{\vec v}(A)) \in \mathcal{I}$, its measure is in $\{0, 1\}$, so the same is true for $A$.

\medskip
This constant $\alpha$ must be $\int f d\mu$, since by definition the expectation has the same integral as $f$ over the whole space, or in a formula
\[ \alpha = \int \alpha d\mu = \int \expectation_\mu(f | \mathcal{I}) d\mu = \int f d\mu. \vspace*{-7mm}\]
\end{proof}

From now on, we concentrate on the hypercubes $H_n = [-n, n]^d$.

\begin{lemma}
\label{lem:DensCharacterization}
Let $X$ be a $\Z^d$-subshift and $f : X \to \R_{+}$ be a simple continuous function. Then there exists $x \in X$ with $\delta_f(x) = \delta_f(X) = \min_{\mu \in \mathcal{M}} \int f(x) d\mu(x)$.
\end{lemma}

\begin{proof}
Denote $\delta_f(X) = \alpha$. Let the radius of $f$ be $R$.
For any $x \in X$ such that $\delta_f(x)$ exists, we can construct a shift-invariant measure $\mu$ such that $\int f(x)d\mu(x) \leq \delta_f(x)$. Namely, let $\mu_n = \frac{1}{(2n+1)^d} \sum_{\vec v \in [-n, n]^d} \delta_{\sigma_{\vec v}(x)}$ where $\delta_y$ denotes the Dirac measure on a configuration $y$ (giving probability $1$ to $y$). The measure $\mu_n$ is just the convex combination of the Dirac measures, i.e.\ picks a random shift of $x$ by some $\vec v \in [-n,n]^d$. Thus the measure $\int f(x)d\mu_n(x)$ is the sum of $f$-values of shifts of $x$ in the hypercube $[-n, n]^d$. Note that this sum only depends on the contents of $x$ in $[-n-R, n+R]^d$.

Taking a weak-$*$ limit point $\mu$ of these measures along the sequence of $n$ reaching the $\liminf$ defining $\alpha$ (using weak-$*$ compactness of the set of measures), we have $\int f(x) d\mu(x) = \alpha$. Furthermore, it is easy to see that $\mu(C) = \mu(\sigma_{\vec u}(C))$ for any clopen set $C$ and any $\vec u \in \Z^d$. Namely, $|\mu_n(C) - \mu_n(\sigma_{\vec u}(C))| \rightarrow 0$ because $\mu_n(C)$ counts the sum of $f$-values of shifts of $x$ in the hypercube $[-n, n]^d$, and $\mu_n(\sigma_{\vec u}(C))$ the same over $[-n, n]^d + \vec u$, and these counts differ by at most $O(n^{d-1})$, so the difference disappears in the normalization by $\frac{1}{(2n+1)^d}$ as $n \rightarrow \infty$.

The above shows that $\inf_{\mu \in \mathcal{M}} \int f(x)d\mu(x) \leq \alpha$. We now observe that $\int f(x)d\mu(x)$ is continuous in $\mu$ (by the definition of the topology), so because $\mathcal{M}$ is weak-$*$ compact, the minimum is reached by some measure $\mu$, showing $\min_{\mu \in \mathcal{M}} \int f(x)d\mu(x) \leq \alpha$.

Next, let $\mu$ be any shift-invariant measure. Recall that $\mathcal{M}$ is spanned by its ergodic measures, meaning that $\mu$ can be written as a (continuous) convex combination of ergodic measures $\nu$, i.e.\ measures where every invariant Borel set has measure $0$ or $1$. In a formula, $\mu(A) = \int \mu_y(A) d\nu(y)$ for all Borel sets $A$, for some probability measure $\nu$ on the set of ergodic measures $\mu_y$. See \cite[Theorem~14.10]{Ge11}. In particular, some ergodic measure $\mu_y$ must satisfy $\int f(x)d\mu_y(x) \leq \alpha$ as well.

The pointwise ergodic theorem for $\Z^d$-actions says that for an ergodic measure $\mu_y$, for $\mu_y$-almost all configurations $x \in X$ we have $\delta_f(x) = \int f(x)d\mu_y(x)$, showing $\alpha \leq \int f(x)d\mu_y(x)$.
Then $\alpha = \int f(x)d\mu_y(x)$, and at least one $x \in X$ satisfies $\delta_f(x) = \alpha$, as desired.
\end{proof}

When $X$ is a subshift of finite type with the uniform filling property, the infimum is reached by some configuration by a simple argument that avoids the use of ergodic theory (or even measures). Namely, we can reach it by gluing together squares with small integral as in \cite{El98}.

Before getting to our actual method in the next section, we prove in the remainder of this section that just from general principles, the minimum density of any symbol can be computed for most subshifts in coding theory. Of course, these methods require being able to efficiently calculate minimal densities in very large finite graphs.

\begin{lemma}
\label{lem:LowerSemicomputable}
Let $X \subset A^{\Z^d}$ be an effectively closed subshift. Then the minimum density of any $a \in A$ in $X$ is lower semicomputable uniformly in the effective description of $X$ by forbidden patterns.
\end{lemma}

\begin{proof}
Let $\alpha$ be the minimum density. By the definition of $X$, for each $n$, we can calculate better and better upper approximations to the set of $[-n, n]^d$-patterns in $X$. Namely, if $F_1, F_2, F_3, \ldots$ is a computable sequence of forbidden patterns defining $X$, then for each $k \in \N$, let $L(k)$ be the set of patterns $P \in A^{[-k,k]^d}$ in which no translate of any $F_1, F_2, \ldots, F_k$ occurs, and let $L(n,k) = \{ P|_{[-n,n]^d} \mid P \in L(k)\}$. Then $L(n,k)$ for $k \geq n$ is a computable sequence of sets converging down to the set of $(2n+1)^d$-patterns that actually occur in configurations of $X$ (but note that we have no control over the speed of convergence). For each $n$, define
\[ \delta_n = \frac{\min\{|\{\vec v \in [-n, n]^d \;|\; x_{\vec v} = a\}| \;|\; x \in X\}}{(2n+1)^d}. \]
We can calculate arbitrarily good lower approximations to the numbers $\delta_n$ using the sets $L(n,k)$, since the minimum can only increase as the set becomes smaller.

\medskip
Partitioning a configuration $x \in X$ into disjoint translates of $[-n, n]^d$, we see $\alpha \geq \delta_n$ for all $n$. 
On the other hand, we can easily construct a shift-invariant measure $\mu \in \mathcal{M}_\sigma(X)$ such that $\mu([a]) = \liminf_n \delta_n$ as in the proof of Lemma~\ref{lem:DensCharacterization}. Hence
\[ \alpha = \min_{\mu \in \mathcal{M}} \mu([a]) \leq \liminf_n \delta_n. \]

Now $\delta_n \leq \alpha \leq \liminf_n \delta_n$ implies $\lim \delta_n = \alpha$. Thus, lower approximations of $\delta_n$ over all $n$ provide converging lower approximations to $\alpha$.
\end{proof}

\begin{lemma}
\label{lem:UpperSemicomputable}
Let $X \subset A^{\Z^d}$ be a subshift of finite type. Suppose that one of the following conditions holds:
\begin{itemize}
\item the language of $X$ is computable, $X$ is block-gluing, and a gluing radius $r$ is known, or
\item $X$ has a totally periodic point and the uniform filling property.
\end{itemize}
Then the minimum density of any $a \in A$ in $X$ is uniformly computable in the description of $X$.
\end{lemma}

More precisely, we show that there exists a Turing machine such that given a Turing machine that describes the language of a subshift which is block-gluing and of gluing radius $r$, it enumerates upper and lower approximations to the minimum density (assuming this data is correct). Under the second condition, we show a bit more: given the description of any SFT $X$, if $X$ happens to have a totally periodic point and happens to have the uniform filling property, then the algorithm correctly computes the density in the sense that it eventually enumerates all lower and upper approximations of the density, and even if it is given an SFT \emph{without} these properties, it never outputs an incorrect approximation.

For example, any subshift of finite type with a \emph{safe symbol} (one contained in no forbidden pattern) satisfies both conditions. When $d \leq 2$, block gluing implies\footnote{For a proof that the uniform filling property implies this, see \cite{Li03}. For our stronger claim, glue together a vertical stack of copies of a horizontal strip. Since we are in an SFT, we can use the same gluing between any two of them, and thus we obtain a vertically periodic point. Then we recall that a two-dimensional SFT with a vertically periodic point has a totally periodic point.} the existence of a totally periodic point; in particular in dimension two we can compute the minimum density of any subshift of finite type with the uniform filling property.

\begin{proof}
Let the minimum density be $\alpha$. Suppose the first condition holds. Since the language is computable, we can compute the minimal density $\delta_n$ of a pattern on $[-n, n]^d$. As seen in the proof of Lemma \ref{lem:LowerSemicomputable}, $\alpha \geq \delta_n$, so we can compute arbitrarily good lower approximations to $\alpha$. On the other hand, if $P$ has density $\delta_n$, then by block-gluing with radius $r$, $X$ contains at least one configuration where $\vec v + [-n, n]^d$ contains a shifted copy of $P$ for each $\vec v \in (2n+1+r)\Z^d$, and by Lemma~\ref{lem:DensCharacterization} applied to the action of $(2n+1+r)\Z^d$, there is also such a configuration where the density exists. Upper bounding the density by assuming all of the ``gluing coordinates'' contain the symbol $a$, we obtain
\[ \alpha \leq \frac{(2n + 1)^d \delta_n + ((2n + 1 + r)^d - (2n + 1)^d)}{(2n + 1 + r)^d}. \]

Now note that $(2n + 1 + r)^d - (2n + 1)^d = O(n^{d-1})$ if $r$ is considered constant. It follows that as $n \rightarrow \infty$, the effect of $r$ becomes negligible. As we also have $\delta_n \rightarrow \alpha$ (by the proof of the previous lemma), the upper bounds given by this formula give arbitrarily good upper approximations to $\alpha$. (Note that the algorithm need not estimate how close the upper bounds are to $\alpha$ at each point, all we need is that both the upper and lower approximations converge to $\alpha$. However, the algorithm does need to evaluate the formula above for each $n$, which requires knowledge of $r$.)

Suppose next that the second condition holds. We conclude lower semicomputability from Lemma \ref{lem:LowerSemicomputable} and now need to show upper semicomputability. We show that under the assumptions, densities of periodic points provide converging upper approximations to $\alpha$, which shows its upper semicomputability since the densities of totally periodic points can be easily enumerated in a subshift of finite type.

Let $x \in X$ be any totally periodic point. Then, for any pattern $P \in A^{[-n, n]^d}$ in the language of $X$, we can successively glue $P$ into a periodic set of positions in $x$ using UFP so that the occurences are disjoint. If the lattice of positions where $P$ is glued aligns with the period lattice of $x$ and the positions of the lattice are separated by at least $2n+1+3r$ where $r$ is the maximum of the UFP radius and the size of a forbidden pattern of $X$ (using that $X$ is of finite type), then we can ensure that the resulting configuration is also totally periodic, and a calculation similar to the previous case shows that the densities converge to $\alpha$. (Note that the algorithm does not need to construct periodic points according to this procedure, although this is one possibility; it suffices that these periodic configurations exist.) 
\end{proof}

In the case where $X$ has UFP and has a totally periodic point, the proof shows that the minimum density is approximated arbitrarily well by periodic points. Most subshifts of coding-theoretic origin have these properties. We note, however, that there are subshifts coming from identifying codes, where the minimum density is not reached by any periodic point \cite{Ob03} (the graphs are gridlike, and are constructed so that identifying codes correspond to valid tilings by a Wang tile set).

\section{Density in periodic configurations}

Recall that a \emph{simple cycle} in a graph $G$ is a cycle that does not repeat vertices.

\begin{lemma}
\label{lem:ReachedByCycle}
Let $G = (V,E,\lambda)$ be a finite directed weighted graph. Then the minimum density of a bi-infinite walk on $G$ is reached by a simple cycle.
\end{lemma}

\begin{proof}
By translating, we may assume the minimum weight of an edge is nonnegative.
Let $\alpha$ be the minimum density of a bi-infinite walk.
We first show that $\alpha$ is also the infimum of the densities of all cycles in $G$.
Let $x$ be a bi-infinite walk of density $\alpha$.
For all $n \in \N$, the walk $w = x_{[-n,n]}$ can be cut into cycles plus at most $|V|$ additional edges.
Namely, consider the leftmost node $q \in V$ in $w$.
If it occurs in $w$ at least twice, then we can remove from $w$ a prefix that also ends in $q$, which is a cycle.
If $q$ occurs only once, we remove the leftmost edge of $w$.
Then we continue with the remaining part of $w$.
The second case can occur only once per each node of $V$, so at most $|V|$ times in total.
All in all we have written $w$ as a sequence of at most $k \leq |V|$ cycles $c_i$, with $t \leq |V|$ transitions in between.

\medskip
Now, writing $\lambda(u)$ for the sum of weights of edges of a walk $u$, we have $\lambda(w) \geq \sum_i \lambda(c_i)$, so for the density we get
\[ \frac{\lambda(w)}{|w|} \geq \frac{\sum_i \lambda(c_i)}{|w|} = \left(\frac{\sum_i \lambda(c_i)}{\sum_i |c_i|} \right)\left( \frac{\sum_i |c_i|}{|w|}\right) \geq
\left(\frac{\sum_i \lambda(c_i)}{\sum_i |c_i|} \right)\left( \frac{2n+1-|V|}{2n+1}\right) \]
from which
\[ \frac{\sum_i \lambda(c_i)}{\sum_i |c_i|} \leq \left(\frac{\lambda(w)}{|w|}\right)\left( \frac{2n+1}{2n+1+|V|}\right). \]
The left side is an affine combination of the weights $\lambda(c_i)/|c_i|$, so for each $n$ one of these cycles must have density below the upper bound. Since $\frac{2n+1}{2n+1+|V|} \rightarrow 1$, the claim follows.

\medskip
We now claim that $\alpha$ is also the infimum density of all simple cycles on $G$. Then it will be also the minimum density of such cycles, since there are only finitely many simple cycles in $G$. Since $\alpha$ is the infimum of densities of cycles, it suffices to show that the infimum of densities of simple cycles is smaller than or equal to that of general cycles. For this, let $c$ be any cycle. 
If $c$ is not simple, it repeats a vertex and we can dissect it into two strictly shorter cycles, $c_1$ and $c_2$.
The density of $c$ is an affine combination of the densities of $c_1$ and $c_2$, so one of them is less than or equal to the density of $c$.
We repeat this procedure until a simple cycle remains; its density is at most that of $c$.
\end{proof}

\begin{lemma}
\label{lem:1DSFTdensity}
If $(Z, \sigma)$ is a one-dimensional subshift of finite type, and $f : Z \to \Q_+$ is a simple continuous function, then the minimum density of $f$ in $Z$ is a rational number, exactly computable from the descriptions of $Z$ and $f$.
\end{lemma}

\begin{proof}
We simply perform a recoding argument to reduce to Lemma \ref{lem:ReachedByCycle}. Specifically, any subshift of finite type is conjugate to the set of paths in a graph: take the edges to represent words of length $k+1$ in the SFT, and nodes words of length $k$, where the maximal length of a forbidden word is at most $k+1$.\footnote{This construction corresponds to the \emph{higher block presentation} in symbolic dynamics (see \cite[Chapter~2.3]{LiMa95}).} If $k$ is further larger than the radius of $f$, then we can represent $f$-images by edge weights.

\medskip
By Lemma \ref{lem:ReachedByCycle}, the minimum density of a configuration is reached by a simple cycle, and since there are a finite number of those, they can be enumerated algorithmically.
\end{proof}

Of course, the minimum density of a simple cycle in a weighted graph can be computed much more efficiently than by exhaustive enumeration, for example by Karp's algorithm \cite{Ka78}.

\begin{remark}
Recall that a $d$-dimensional \emph{sofic shift} is a subshift of the form $\pi(X)$, where $X \subset A^{\Z^d}$ is an SFT and $\pi : X \to B^{\Z^d}$ is a relabeling map, defined by $\pi(x)_{\vec n} = \hat\pi(x_{\vec n})$ for all $\vec n \in \Z^d$ for some symbol-to-symbol function $\hat\pi : A \to B$. Lemma \ref{lem:1DSFTdensity} is also true for one-dimensional sofic shifts, as we can simply replace a sofic shift $Z = \pi(X)$ by the SFT $X$, and the function $f$ by $f \circ \pi$.
\end{remark}

We now transfer the computability of densities in one-dimensional subshifts to periodic parts of multidimensional subshifts.

\begin{lemma}
\label{lem:PeriodicComputableDensity}
Let $Y \subset A^{\Z^d}$ be a subshift of finite type, and let $f : Y \to \Q_+$ be a simple continuous function. Let $L \subset \Z^d$ be a rank-$(d-1)$ subgroup. Then the minimal density of an $L$-periodic configuration in $Y$ is a rational number, exactly computable from the descriptions of $Y$, $f$ and $L$.
\end{lemma}

\begin{proof}
Let $\vec v \in \Z^d$ be such that $\hat L = \langle L, \vec v \rangle$ has rank $d$. Let $F$ be a (finite) set of coset representatives for $\hat L$, meaning $\hat L + F = \Z^d$, and each $\vec u \in \Z^d$ has a unique representation as such a sum.

Let  $Y_L$ be the set of $L$-periodic configurations. The minimum density of an $L$-periodic configuration is $\inf_{\mu \in \mathcal{M}_\sigma(Y_L)} \int f(x) d\mu(x)$, where $\mathcal{M}_\sigma(Y_L)$ is the set of shift-invariant measures on $Y_L$ (note that we mean invariance with respect to $\Z^d$-shifts, not just those in $L$). To see this, apply Lemma~\ref{lem:DensCharacterization} applied to $Y_L$, seen as a subshift under the shift-action of $\Z^d$

\medskip
Now, letting $g(x) = \sum_{\vec u \in F} f(\sigma_{\vec u}(x)) / |F|$, we claim that
\[ a = \inf_{\mu \in \mathcal{M}_\sigma(Y_L)} \int f(x) d\mu(x) = \inf_{\mu \in \mathcal{M}_{\hat L}(Y_L)} \int g(x) d\mu(x) = b \]
where $\mathcal{M}_{\hat L}(Y_L)$ denotes the $\hat L$-invariant measures on $Y_L$ (meaning measures invariant under translations $\sigma_{\vec u}$ with $\vec u \in \hat L$). To see this, for $a \geq b$ we note that any measure considered in the infimum on the left (i.e.\ any invariant measure) is also considered by the infimum on the right (since $\Z^d$-invariance implies $\hat L$-invariance) and gives the same integral on both sides. For $a \leq b$ we note that for any measure $\mu$ considered on the right, the convex combination of measures $\sum_{\vec u \in F} \sigma_{\vec u}(\mu) / |F|$ is considered on the left and gives the same integral. Here, the measure $\nu = \sigma_{\vec u}(\mu)$ is defined by $\nu(B) = \mu(\sigma_{-\vec u}(B))$ for $B$ Borel.

Now consider the action of the shift by $\vec v$ on the $L$-periodic configurations $Y_L$ in $Y$. It is easy to see that $\sigma_{\vec v}(Y_L) = Y_L$, and the system $(Y_L, \sigma_{\vec v})$ is topologically conjugate to a subshift of finite type, by the map $\phi : (Y_L, \sigma_{\vec v}) \to (Z, \sigma) \subset ((A^F)^\Z, \sigma)$ defined by $(\phi(x)_i)_{\vec u} = x_{i \vec v + \vec u}$.

\medskip
Now since $Y_L$ is \emph{pointwise} $L$-invariant by definition, we have
\[ \inf_{\mu \in \mathcal{M}_{\hat L}(Y_L)} \int g(x) d\mu(x) = \inf_{\mu \in \mathcal{M}_{\sigma_{\vec v}}(Y_L)} \int g(x) d\mu(x). \]
But as discussed $(Y_L, \sigma_{\vec v})$ is topologically conjugate to some SFT $(Z, \sigma)$ and $g$ corresponds to some~simple continuous function $h$. By Lemma \ref{lem:1DSFTdensity}, this value is rational and can be effectively computed.
\end{proof}

\begin{definition}
The \emph{hexagonal grid} is the undirected graph $(\Z^2, E)$ where $E$ contains all edges of the forms $((x,y), (x+1,y))$ and $((x,y), (x,y+p))$, where $p = (-1)^{x+y}$, for $(x, y) \in \Z^2$.
\end{definition}

In words, we include all horizontal edges between adjacent coordinates, and we include a vertical edge between $(x, y)$ and $(x, y+1)$ if and only if $x+y$ is even. It is easy to see that this grid is isomorphic (as a graph) with the hexagonal grid as depicted in Figure~\ref{fig:sharpSAT} and Figure~\ref{fig:NewRecord}. 

The hexagonal grid is a gridlike graph.
Namely, translations by the group $T = \langle (0, 2), (1, 1) \rangle \simeq \Z^2$ are graph automorphisms, and divide $\Z^2$ into two orbits based on the parity of the coordinates. The set $X_{\hex} \subset \{0,1\}^{\Z^2}$ of identifying codes on the graph $(\Z^2, E_{\hex})$ then corresponds to a $\Z^2$-subshift of finite type $Y_{\hex}$ over the alphabet $\{0,1\}^2$. The forbidden patterns of $Y_{\hex}$ are induced by subsets of $\Z^2$ that cannot all contain a 0 in a valid identifying code, either because some node would have no codewords in its neighborhood, or because two nodes would have identical identifying sets.

\begin{figure}[!h]
\vspace*{2mm}
\input{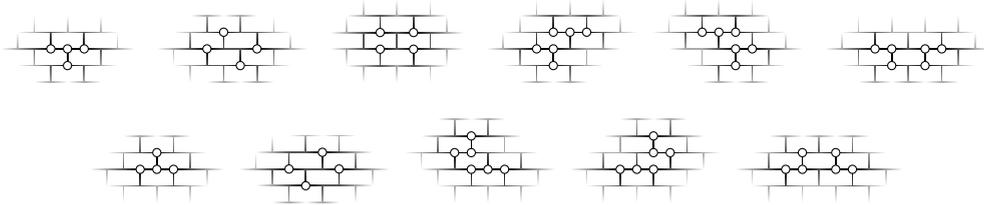}
\caption{The minimal set $F_{\hex}$ of forbidden patterns for identifying codes on $E_{\hex}$. The circles denote non-codewords.}
\label{fig:hexforbs}
\end{figure}

More explicitly, first we observe that we can pick a finite set $F_{\hex}$ of sets, which define $X_{\hex} \subset \{0,1\}^{\Z^2}$ when for each $T$-translate $F'$ of each $F \in F_{\hex}$ we forbid the all-zero pattern $p \in \{0\}^F$. Such a family of sets is shown in Figure \ref{fig:hexforbs}. Picking the leftmost (and in case of ties the bottommost) cell of a subset of $\Z^2$, as its ``anchor'', we have
\[ F_{\hex} \!=\! \{\{(0,0), (1,0), (1,-1), (2,0)\}, \{(0,0), (1,1), (2,-1), (3,0)\}, \{(0,0), (0,1), (2,0), (2,1)\}\! ... \}\]
(with eight sets omitted, but copiable from the figure). Now, to define $Y_\hex$, we induce a linear isomorphism $\ell : \Z^2 \to T$ by $\ell((1,0)) = (1,1); \ell((0,1)) = (0, 2)$, pick the set of coset representatives $\{(0,0), (1,0)\}$ for $T$, and then define $f : \{0,1\}^{\Z^2} \to (\{0,1\}^2)^{\Z^2}$ by $f(x)_{\vec v} = (x_{\ell(\vec v)}, x_{\ell(\vec v) + (1,0)})$. Now $Y_\hex = f(X_{\hex})$, and each pattern in $F_{\hex}$ corresponds to a finite set of forbidden patterns defining $Y_{\hex}$.

\medskip
From Lemma~\ref{lem:PeriodicComputableDensity} we now obtain the following:

\begin{lemma}
\label{lem:CanCompute}
Given $\vec v \in T$, the minimal density of a $\vec v$-periodic identifying code on $E_{\hex}$ is a rational number exactly computable from $\vec v$.
\end{lemma}

Most of the codes in the literature are found by looking at two periods, i.e.\ a period lattice of full rank. Both methods lead to the same upper bounds on the minimal density, since as shown in the proofs of Lemma~\ref{lem:PeriodicComputableDensity} and Lemma~\ref{lem:ReachedByCycle}, the minimal density of a configuration with a fixed period lattice of rank $d-1$ is actually reached by a configuration with period lattice of rank $d$.
However, heuristically our method should lead to a more efficient (or at least essentially different) algorithm, as we do not need to guess in advance how large the period complementary to $\vec v$ needs to be. And indeed, as reported in Theorem~\ref{thm:NewRecord}, we were able to improve on the smallest known density of an identifying code this way.

\section{Implementation}

By Lemma~\ref{lem:CanCompute}, we can compute the exact minimal density of any $\vec v$-periodic identifying code on $E_{\hex}$. Specifically, this reduces to computing the minimum density of a cycle in a certain weighted directed graph $G = (V,E,\lambda)$ with integer weights.

We describe our construction of the graph $G$.
Let $(x,y) = \vec v$ and assume $y > 0$, $x+y$ even and $x \geq 0$.
For $j \in \Z$, let $s_j = \lceil j x/y \rceil$, and define the \emph{border} as $B = \{(s_j,j), (s_j+1,j) \mid j \in \Z\}$.
The nodes of our graph are sets $C$ of translates of the sets $F_{\hex}$ that are invariant under translation by $\vec v$, meaning that $\{ S + \vec v \mid S \in C \} = C$, and such that each $S \in C$ contains both a node $(i,j)$ with $i \leq s_j+1$ and a node $(i,j)$ with $i \geq s_j$; we call this the \emph{border condition}.
Note that the sets $C$ are technically infinite (if nonempty), but the periodicity condition ensures they have a finite encoding (i.e.\ can be described with a finite amount of information).

\medskip
In our graph, we include an edge $(C,D) \in E$ whenever there exists a $\vec v$-periodic pattern $P \in \{0,1\}^B$ with
\[
  D = \{ S - (2,0) \mid S \in C \cup C_B, \nexists \vec n \in S: P_{\vec n} = 1 \}
\]
where $C_B$ contains all translates of sets in $F_{\hex}$ that intersect $B$ and are contained within the half-plane $\{ (i,j) \in \Z^2 \mid s_j \leq i\}$.
The weight of such an edge is $\sum P$, or the number of 1-symbols in $P$.
The idea is that a bi-infinite walk in $G$ corresponds to a $\vec v$-periodic configuration specified one (slanted) vertical strip -- a translate of $B$ -- at a time.
An edge $(C,D) \in E$ corresponds to an adjacent pair of such strips.
The nodes maintain a collection of forbidden patterns that can be ``dealt with'' as soon as a 1-symbol, representing a codeword, is known to occur in them.
If we have not been able to get rid of a forbidden pattern, it is shifted to the left.
Patterns are forbidden from drifting to the left of the border $B$, as they would then become impossible to remove. (Formally, one can still define $D$ by the same formula, but $D$ will not satisfy the border condition, so one does not obtain an edge in the graph.)

\medskip
Again, since the states are $\vec v$-periodic and contain only translates of sets in $F_{\hex}$ that satisfy the border condition, the number of possible states is finite and the edges can be easily computed.
We construct our graphs by starting from the empty node $C_0 = \emptyset$ and including every node reachable from it.
The resulting subgraph is strongly connected (as $C_0$ is reachable from every node by all-1 edges) and contains exactly the same labels of bi-infinite walks as the full graph.

Karp's algorithm \cite{Ka78} can be used to calculate the minimum density of a cycle in the graph $(G,E,\lambda)$. We chose this algorithm because it is easy to implement and readily parallelizable.
We did not use any special tricks in the construction of the graph, apart from basic memory optimization and parallelization. Karp's algorithm was the bottleneck in all computations, and we used some tricks to be able to go further in the calculations. We explain these tricks in this section.

\medskip
First, let us recall the basic Karp's algorithm for a weighted graph $G = (V,E,\lambda)$. The algorithm only works when the graph has an \emph{initial node} $s$ from which every other node is reachable. This can be guaranteed by passing to strongly connected components, but in our application $G$ is always strongly connected. Now, for each node $v$ and length $0 \leq k \leq n$, compute the minimal weight $F_k(v)$ of path from $s$ to $v$. This can be done using the recurrence $F_k(v) = \min_{(w,v) \in E} (F_{k-1}(w) + \lambda(w,v))$. The minimum density of a cycle is then
\begin{equation}
  \label{eq:karp}
  \alpha = \min_{v \in V} \max_{0 \leq k < n} \frac{F_n(v) - F_k(v)}{n - k}.
\end{equation}
This can be computed in space $O(n^2)$ where $n = |V|$ by storing the weights $F_k(v)$ in an $n \times (n+1)$ array, or in space $O(n)$ by using a handful of length-$n$ arrays: first, compute the weights $F_k(v)$ for $k = 0, 1, \ldots, n$ using two arrays, then store the weights $F_n(v)$ into a third array $A$, and finally recompute the weights $F_k(v)$ using two arrays while simultaneously computing \eqref{eq:karp} using the values of the array $A$ and a fourth array that stores the running maximum of the fractions $(F_n(v) - F_k(v))/(n-k)$ for each $v \in V$. Both variants run in $O(n^2)$ time, but the latter takes approximately twice as long.

\medskip
In order to actually compute a cycle of minimum density, one takes a node $v^* \in V$ that attains the minimum in \eqref{eq:karp} and computes a length-$n$ path $P$ from $s$ to $v^*$ of minimum weight.
It necessarily contains a minimum density cycle as a contiguous subpath, and it is easy to locate one knowing the value of $\alpha$.
The path $P$ can be computed as part of the $O(n^2)$ space algorithm by storing, in addition to each value $F_k(v)$, a node $w$ with $(w,v) \in E$ that minimizes $F_{k-1}(w) + \lambda(w,v)$, and following these nodes backward from $v^*$ to $s$.

\medskip
However, $O(n^2)$ space proved to be too high for the code of density 53/126 in Figure \ref{fig:NewRecord}; we could prove its existence using the $O(n)$ variant, but not extract a concrete code, as the cluster ran out of memory.
For this purpose we implemented a third variant that runs in space $O(n^{3/2})$ and takes about three times as long as the $O(n^2)$ variant:
\begin{enumerate}
\item Fix a starting node $s \in V$.
\item Pick $m = \lceil \sqrt{n} \rceil$ numbers $0 = n(1) < n(2) < \cdots < n(m) = n$ with $n(i+1) \leq n(i) + m$ for each $i$.
\item For each $k = 0, \ldots, n$, compute the weights $F_k(v)$ for $v \in V$ using two length-$n$ arrays. If $k = n(i)$ for some $i$, then store these weights into one row of an $m \times n$ array $A$.
\item For each $k = 0, \ldots, n$, recompute the weights $F_k(v)$ and simultaneously compute $\alpha$ from \eqref{eq:karp} and the minimizing vertex $v^*$.
\item Set $P$ as the path containing only $v^*$.
\item For each $i = m, m-1, \ldots, 2$, compute an $m \times n$ array $B$ containing, for each $n(i-1) < k \leq n(i)$ and $v \in V$, the vertex $w$ with $(w,v) \in E$ that minimizes $F_{k-1}(w) + \lambda(w,v)$.
  The values $F_{n(i-1)}(v)$ stored in $A$ are used as a starting point.
  From this data, extract a length-$m$ path $Q$ from some $v \in V$ to the beginning of $P$ that minimizes $F_{n(i-1)}(v) + \lambda(Q)$.
  Prepend it to $P$.
\end{enumerate}
At the end of the algorithm, $P$ is a length-$n$ path from $s$ to $v^*$ of weight $F_n(v)$.

\subsection*{Acknowledgments}

We thank Ville Junnila and Tuomo Lehtil\"a  for suggesting problems to apply our method on. We thank Tim Austin for suggesting a reference for the pointwise ergodic theorems specific to the abelian setting. The computer resources of the Finnish IT Center for Science (CSC) and the FGCI project (Finland) are acknowledged.


\end{document}